\documentclass[11pt]{article}
\usepackage{amssymb}
\usepackage{amsmath}
\usepackage{amsthm}
\usepackage{latexsym}
\usepackage{amsfonts}
\usepackage{graphicx}
\usepackage{graphics}

\newcommand{\dis}{\displaystyle}
\theoremstyle{definition}
\newtheorem{Def}{Definition}[section]
\theoremstyle{plain}
\newtheorem{thm}[Def]{Theorem}  
\newtheorem{prop}[Def]{Proposition}

\theoremstyle{definition}

\newtheorem*{Proof}{Proof}

\newcommand{\el}{\ell}

\newcommand{\ra}{\;\rightarrow\;}

\newcommand{\al}{\alpha}

\newcommand{\ga}{\gamma }
\newcommand{\Ga} {{\varGamma}}

\newcommand{\de}{\delta }
\newcommand{\OO} {{\varOmega}}

\newcommand{\C}{\mathbb{C}}

\newcommand{\N}{\mathbb{N}}

\newcommand{\oB}{\overline{B}}

\newcommand{\cu}{{\cal{U}}}

\newcommand{\ce}{{\cal{E}}}

\newcommand{\ld}{\ldots}

\newcommand{\sm}{\smallsetminus}

\newcommand{\qb}{$\quad\blacksquare$}

\begin{document}
\title{Domains of Holomorphy}
\author{V. Nestoridis}
\date{}
\maketitle\vspace*{-1cm}
\begin{center}
Dedicated to my friend Professor Paul M. Gauthier
\end{center}
\bigskip
\begin{abstract}
We give a simple proof that the notions of Domain of Holomorphy and Weak Domain of Holomorphy are equivalent. This proof is based on a combination of Baire's Category Theorey and Montel's Theorem. We also obtain generalizations by demanding that the non-extentable functions belong to a particular class of functions $X=X(\OO)\subset H(\OO)$. We show that the set of non - extendable functions not only contains a $G_{\delta}$ - dense subset of $X(\OO)$, but it is itself a $G_{\delta}$ - dense set. We give an example of a domain in $\C$ which is a $H(\OO)$-domain of holomorphy but not a $A(\OO)$-domain of holomorphy.
\end{abstract}

{\em AMS classification number}\,: Primary 32T05, Secondary 30D45.\vspace*{0.2cm} \\
{\em Key words and phrases}\,: extendability, Domain of holomorphy, Weak Domain of holomorphy, Baire's theorem, generic property, Montel's Theorem, Analytic continuation.
\section{Introduction}\label{sec1}
\noindent

It is well - known that the notions of domain of holomorphy and of weak - domain of holomorphy are equivalent. The original proof is constructive, technical and by no means elementary (\cite{Range}). A simpler proof was obtained in \cite{PP} by a combination of Baire's theorem and a theorem of Banach. Furthermore in \cite{PP} it was proven that the set of non - extendable functions contains a $G_{\delta}$ - dense subset of $H(\OO)$ of holomorphic functions on a domain $\OO$, or more generally in a space $X(\OO) \subset H(\OO)$ satisfying some assumptions. The use of Banach's theorem does not allow to conclude that the set of non - extendable functions is itself a $G_{\delta}$ - dense subset of $X(\OO)$; in fact the sets appearing in Banach's theorem can be very high in Borel's hierarchy (\cite{S.R}, \cite{LD.SG}).

In the present paper we replace Banach's theorem by Montel's theorem and combining it with Baire's theorem we obtain a new, very simple proof of complex analytic nature of a slightly stronger result; that is, we prove that the set of non - extendable functions is itself $G_{\delta}$ - dense and not only that it contains a $G_{\delta}$ - dense set.

At present, $\OO$ is a domain in the finite dimensional space $\mathbb{C}^d$. In future papers we will discuss the infinite dimensional care. We mention that in \cite{Ma} Montel's theorem was used to treat the case $X(\OO) = H(\OO)$, where $\OO$ is a domain in a separable Banach space, even infinite dimensional. For finite dimensional holomorphy Montel's theorem towards generic results has also been used in the works of Paul M. Gauthier; see for instance \cite{Gauthier 2}.

Some generic results for particular choices of $X(\OO)$ have already been obtained in \cite{BNP}, \cite{BNPP}, \cite{Georgakopoulos}.

It is well-known that in $\C$ every domain is a $H(\OO)$-domain of holomorphy. We give an example of a domain in $\C$ which is a $H(\OO)$-domain of holomorphy but it is not a $A(\OO)$-domain of holomorphy, where $A(\OO)=\{f:\overline{\OO}\ra\C$, continuous on $\overline{\OO}$ and holomorphic in $\OO\}$. Such an example is the domain $\OO=\{z\in\C:|z|<1,\;z\notin[0,1/2]\}$. Our proof (in several variables) implies in particular that if $\OO\subset\C^d$ is not a $X(\OO)$-domain of holomorphy, then there exist two balls $B_1,B_2$, $B_2\cap\OO\neq\emptyset$, $B_2\cap\OO^c\neq\emptyset$, $B_1\subset\overline{B}_1\subset B_2\cap\OO$ such that every $f\in X(\OO)$ restricted to $B_1$ admits a bounded holomorphic extension in $B_2$.

Possible particular choices of $X(\OO)$ are the spaces $A^p(\OO)$, $H^\infty_p(\OO)$, $p\in\{\infty\}\cup\{0,1,2,\ld\}$, the Bergman spaces $B_p(\OO)$, $0<p<+\infty$, the Hardy spaces and variations or combinations of them in one or several complex variables. Here $A^p(\OO)=\{f:\OO\ra\C$ holomorphic, such that the derivative $f^{(\el)}$ extends continuously in $\overline{\OO}$ for all $\el\in\{0,1,2,\ld\}$, $\el\le p\}$ and $H^\infty_p(\OO)=\{f:\OO\ra\C$ holomorphic such that each derivative $f^{(\el)}$ is bounded on $\OO$ for all $\el\in\{0,1,2,\ld\}$, $\el\le p\}$. In the case of $X=A^p(\OO)$, under some assumptions, there is a relation with the $p$-continuous analytic capacity $\al_p$ introduced in \cite{BNPP}. Similarly for the case $X=H^\infty_p(\OO)$, I believe that it is possible to define an analogous notion of $p$-analytic capacity $\ga_p$ which relates to our situation. It suffices to replace the space $A^p(\OO)$ by the space $H^\infty_p(\OO)$ in the definition of $\al_p$.
\section{Preliminaries}\label{sec2}
\noindent

In \cite{Range} we find the following definition.
\begin{Def}\label{Def2.1}
Let $\OO\subset\C^d$ be open and connected and $f:\OO\ra\C$ be a holomorphic function. Then, $f$ is called extendable if there exist an open and connected set $U\subset\C^d$ with $U\cap\OO\neq\emptyset$ and $U\cap\OO^c\neq\emptyset$, a holomorphic function $F:U\ra\C$ and a component $V$ of $U\cap\OO$ such that $F|_V=f|_V$. Otherwise $f$ is called non-extendable.
\end{Def}

The following definition comes from the theory of analytic continuation mainly in one complex variable.
\begin{Def}\label{Def2.2}
Let $\OO\subset\C^d$ be open and connected and $f:\OO\ra\C$ be a holomorphic function. Then, $f$ is called extendable in the sense of Riemann domains, if there exist two open Eucledian balls $B_1,B_2\subset\C^d$, with
\[
B_1\subset\oB_1\subset B_2\cap\OO, \ \ B_2\cap\OO\neq\emptyset, \ \ B_2\cap\OO^c\neq\emptyset,
\]
and a bounded holomorphic function $F:B_2\ra\C$ such that $F|_{B_1}=f|_{B_1}$. Otherwise the function $f$ is called non-extendable in the sense of Riemann domains.
\end{Def}
\begin{prop}\label{prop2.3}
Definitions \ref{Def2.1} and \ref{Def2.2} are equivalent.
\end{prop}
\begin{Proof}
Suppose that $f$ satisfies Definition \ref{Def2.2}. Set $U=B_2$ and let $V$ the component of $U\cap\OO$ containing $B_1$. Since $f|_{B_1}=F|_{B_1}$, by analytic continuation it follows $f|_V=F|_V$. Thus, Definition \ref{Def2.1} is also satisfied.

Conversely, suppose that $f$ satisfies Definition \ref{Def2.1}. We claim that $\overline{V}\cap\partial\OO\neq\emptyset$. Let $z_1\in V\cap\OO\subset U\cap\OO$ and $z_2\in U\cap\OO^c$. Since $U\subset\C^d$ is open and connected, there is a polygonal line $\Ga$ in $U$ joining $z_1$ and $z_2$. Since $z_1\in V$ and $z_2\notin V$ it follows that $\Ga$ meets $\partial V$. Let $w\in\Ga\cap\partial V$. Since $\Ga\subset U$ it follows easily that $w\in\partial\OO$. Thus, $w\in\partial V\cap\partial\OO\subset\overline{V}\cap\partial\OO$ and the claim is proven.

Since $w\in U$ there is $r>0$, so that $B(w,r)\subset U$, where $B(w,r)$ denotes the Euclidian ball centered at $w$ with radius $r$. We set $B_2=B\Big(w,\dfrac{r}{2}\Big)$. Then $F$ is holomorphic and bounded on $B_2$. Since $w\in\overline{V}$ it follows that $B_2\cap V\neq\emptyset$ and let $z\in B_2\cap V$. Then there is $\de>0$ so that $B(z,\de)\subset B_2\cap V$. We set $B_1=B\Big(z,\dfrac{\de}{2}\Big)$. Thus, $B_1\subset\overline{B}_1\subset B_2\cap V\subset B_2\cap\OO$, $F|_{B_1}=f|_{B_1}$ and $F$ is holomorphic and bounded on $B_2$. Thus, $f$ satisfies Definition \ref{Def2.2}.  \qb
\end{Proof}
\noindent
{\bf Remark}. Often a function which is non extendable is called nowhere extendable, and a function which is extendable is called somewhere extendable.
\section{The result}\label{sec3}
\noindent

Let $\OO\subset \mathbb{C}^d$ be open and connected and $X=X(\OO)$ let be a set of holomorphic functions $f:\OO\ra\C$; that is, $X\subset H(\OO)$.
\begin{Def}\label{Def3.1}
The open connected set $\OO\subset\C^d$ is called an $X$-domain of holomorphy if there exists $f\in X$ which is non-extendable.
\end{Def}
\begin{Def}\label{Def3.2}
The open connected set $\OO\subset\C^d$ is called weak $X$-domain of holomorphy if for every pair of open Eucledian balls $B_1,B_2$ with $B_2\cap\OO\neq\emptyset$, $B_2\cap\OO^c\neq\emptyset$, $B_1\subset\overline{B}_1\subset B_2\cap\OO$ there exists a function $f_{B_1,B_2}\in X$ such that the restriction of $f_{B_1,B_2}$ on $B_1$ does not have any bounded holomorphic extension on $B_2$.
\end{Def}
\begin{thm}\label{thm3.3}
We suppose that $X=X(\OO)\subset H(\OO)$ is a topological vector space endowed with the usual operations $+,\cdot$ and that its topology is induced by a complete metric. We also suppose that the convergence $f_n\ra f$ in $X$ implies the pointwise convergence $f_n(z)\ra f(z)$ for all $z\in\OO$. Then definitions \ref{Def3.1} and \ref{Def3.2} are equivalent. If the above assumptions hold and $\OO$ satisfies definitions \ref{Def3.1} and \ref{Def3.2}, then the set $\{f\in X: f$ is non extendable$\}$ is a dense and $G_\de$ subset of $X$.
\end{thm}
\begin{Proof}
It is obvious that Definition \ref{Def3.1} implies Definition \ref{Def3.2}; it suffices to set $f_{B_1,B_2}=f$. In order to prove the rest it suffices to assume that $\OO$ satisfies Definition \ref{Def3.2} and prove that the set $A=\{f\in X:f$ is non-extendable$\}$ is dense and $G_\de$ in $X$. Equivalently, it suffices to show that $A^c$ is a denumerable union of closed sets in $X$ with empty interiors.

Consider the set of couples $(B_1,B_2)$ of open Eucledian balls such that $B_2\cap\OO\neq\emptyset$, $B_2\cap\OO^c\neq\emptyset$, $B_1\subset\overline{B}_1\subset B_2\cap\OO$ and the centers of $B_1,B_2$ belong to $Q^{2d}$ and the radii of $B_1,B_2$ belong to $(0,+\infty)\cap Q$, where the set $Q$ denotes the set of rational numbers. This set $Y$ is denumerable.

It is easy to see that
\[
A^c=\bigcup_{(B_1,B_2)\in Y}\bigcup_{M\in\{1,2,\ld\}}\{f\in X:\;\exists\; F\;\text{holomorphic on}\; B_2\; \text{and bounded by $M$ so that}\;
\]
$F|_{B_1}=f|_{B_1}\}$.

Thus, $A^c$ is a denumerable union of sets of the form $T(B_1,B_2,M)=\{f\in X:\exists\; F$ holomorphic on $B_2$ bounded by $M$ such that $F|_{B_1}=f|_{B_1}\}$. By Baire's Theorem it suffices to prove that for any fixed choice of $B_1,B_2$ and $M$ the set $T(B_1,B_2,M)$ is closed in $X$ and its interior is empty.

Suppose $f_n\in T(B_1,B_2,M)$ and $f_n\ra f$ in the topology of $X$, where $f\in X$. For each $n\in\{1,2,\ld\}$ there exists a holomorphic function $F_n$ on $B_2$ bounded by $M$ such that $F_n|_{B_1}=f_n|_{B_1}$. By Montel's theorem \cite{Gauthier 1} a subsequence $F_{k_n}$ of $F_n$ converses uniformly on compact subsets of $B_2$ towards a function $F$ holomorphic on $B_2$ and bounded by $M$.

Since the convergence $f_n\ra X$ in the topology of $X$ implies pointwise converge in $\OO$ by assumption, it follows that $f|_{B_1}=\dis\lim_n f_n|_{B_1}=\dis\lim_n f_{k_n}|_{B_1}=\dis\lim_nF_{k_n}|_{B_1}=F|_{B_1}$, where $f\in X$ and $F$ is holomorphic in $B_2$ and bounded by $M$. Thus, $f\in T(B_1,B_2,M)$. This proves that $T(B_1,B_2,M)$ is closed in $X$.

Finally, we shall show that the interior of $T(B_1,B_2,M)$ in $X$ is void. Assume that $T(B_1,B_2,M)^0\ni f$ to arrive at a contradiction. By assumption there exists a function $f_{B_1,B_2}\in X$ such that its restriction on $B_1$ does not have any bounded holomorphic extension on $B_2$. Since $f+\dfrac{1}{n}f_{B_1,B_2}\ra f$ in the topology of $X$ and $f$ is in the interior of $T(B_1,B_2,M)$ it follows that for some $n\in\{1,2,\ld\}$ the function $f+\dfrac{1}{n}f_{B_1,B_2}$ belongs to $T(B_1,B_2,M)$. The same holds for the function $f$. Thus, both functions restricted to $B_1$, admit holomorphic extensions on $B_2$ bounded by $M$. Thus, their difference $\dfrac{1}{n}f_{B_1,B_2}$ restricted to $B_1$ admits a holomorphic extension on $B_2$ bounded by $2M$. It follows that the function $f_{B_1,B_2}$ restricted to $B_1$ admits a holomorphic extension on $B_2$ bounded by $2nM$. This contradicts the fact that $f_{B_1,B_2}|_{B_1}$ does not admit any bounded holomorphic extension on $B_2$. Thus, $T(B_1,B_2,M)^0=\emptyset$ and the proof is complete. \qb
\end{Proof}
\section{Remarks}\label{sec4}
\noindent

If $X=H(\OO)$ we have the notions of Domain of holomorphy and Weak Domain of holomorphy. Their equivalence in this case is well known \cite{Range}, \cite{PP}. Our proof is very simple.

It is also known that for $d=1$ every domain $\OO\subset\C$ is a domain of holomorphy. It is very simple to see that it is a weak domain of holomorphy; it suffices to consider the function $\dfrac{1}{z-\zeta}$ for each $\zeta\in\partial\OO$ \cite{Range}, \cite{PP}.

For $d=1$ the domain $\OO=\{z\in\C:|z|<1,\;z\notin\Big[0,\dfrac{1}{2}\Big]\Big\}$ is an $H(\OO)$-domain of holomorphy but not a $A(\OO)$-domain of holomorphy, where $A(\OO)=\{f:\overline{\OO}\ra\C$, continuous an $\overline{\OO}$ and holomorphic in $\OO\}$. This follows easily by Morera's theorem, because if $\ce$ is a straight line and $\cu\subset\C$ an open set, every function continuous on $\cu$ and holomorphic in $\cu-\ce$ is holomorphic in $\cu$ \cite{Rudin}.

Examples of sets $X=X(\OO)$ satisfying the assumptions of Theorem \ref{thm3.3} (and therefore, its conclusion, as well) are the following in one complex variable.

$A(\OO)$, $A^p(\OO)=\{f\in H(\OO):f^{(\el)}$ extends continuously on $\overline{\OO}$ for all $\el\in\N$, $\el\le p\}$, $p\in\{\infty\}\cup\{0,1,2,\ld\}$ $H^\infty_p(\OO)=\{f\in H(\OO):f^{(\el)}$ is bounded in $\OO$ for all $\el\in\N$, $\el\in p\}$ and mixtures of them, where $N=\{0,1,2,\ld\}$. Also the Bergman and Hardy spaces and variants of them.

Let $L\subset\C$ be a compact set satisfying some assumptions. We set $\OO=\C-L$. Then the fact that $\OO$ is an $A^p(\OO)$-domain of holomorphy or not, relates to the $p$-continuous analytic capacity $\al_p(L)$ \cite{BNPP}. The analogous question for $H^{\infty}_p(\OO)$ should relate to the $p$-analytic capacity $\ga_p(L)$ which, I believe, can be defined in a similar way with the definition of $\al_p(L)$; it suffices to replace the space $A^p(\widehat{\C}\sm L)$ by the space $H^\infty_p(\widehat{\C}\sm L)$.

In several variables one can examine the situation for the analogous spaces $A^p(\OO)$, $H^\infty_p(\OO)$, Bergman and Hardy spaces and variants of them.\medskip\\
\noindent
{\bf Acknowledgment.} I would like to thank R. Aron, B. Braga, S. Charpentier, P. M. Gauthier, T. Hatziafratis, M. Maestro, P. Pflug and A. Siskakis for helpful communications.

V. Nestoridis, \\
National and Kapodistrian University of Athens \\
Department of Mathematics\\
Panepistemiopolis \\
157 84 Athens\\
Greece\\
e-mail:vnestor@math.uoa.gr

\end{document}